\documentclass[12pt,a4paper]{article}

\usepackage[utf8]{inputenc}
\usepackage[english]{babel}
\usepackage{amsthm}
\usepackage{amsfonts}
\usepackage{amsmath}
\usepackage{graphicx}

\newcommand{\re}{\mbox{\rm Re}\,}
\newcommand{\im}{\mbox{\rm Im}\,}
\newcommand{\Ai}{\mbox{\rm Ai}\,}
\newcommand{\Bi}{\mbox{\rm Bi}\,}

\newtheorem{Theorem}{Theorem}

\title{Eigenvalue Dynamics of a $\mathcal{PT}$--symmetric Sturm-Liouville Operator. Criteria of the Similarity
 to a Self-adjoint or Normal Operator.}
\author{Shkalikov A.A., Tumanov S.N.}
\date{25 July 2017}

\begin{document}

\maketitle
%\begin{abstract}
% ...
%\end{abstract}
\frenchspacing

{\bf Introduction.} The goal of the paper is to investigate the dynamics of the eigenvalues
of the Sturm-Liouville operator
\begin{equation}\label{1}
T(\varepsilon)y=-\varepsilon^{-1}y''+p(x)y
\end{equation}
on the finite interval  $[-1,1]$  as the parameter  $\varepsilon >0$ changes.  For simplicity we consider
the Dirichlet boundary conditions
\begin{equation}\label{2}
y(-1)=y(1) =0.
\end{equation}
It is assumed in the sequel that the potential $p$ is summable and $\mathcal{PT}$-symmetric, i.e. $p(x)= \overline{p(-x)}$. It is easily seen that  this condition guarantees the symmetry with respect to the real
axis of the spectrum of the operator $T(\varepsilon)$ defined by differential expression  \eqref{1} and boundary conditions  \eqref{2}.

Plenty of papers are dedicated to the study of  $\mathcal{PT}$-symmetric operators especially
in the physical literature. We point out the papers  \cite{BB1}-\cite{EG} which we acquainted with, and which have stimulated our investigation.  More details and references  can be found in the review of Dorey, Dunning and Tateo
\cite{DDT2}. The papers which are close to the topic of our paper  are dedicated, basically, to the proof of the reality of the spectrum of the 1-d Schrodinger operator on the whole line with some concrete  potentials, in particular,
with cubic or polynomial potentials for particular  parameter values.

In this work we pose the question not only about the realness of the spectrum, but also on the similarity of the operator
to a self-adjoint or normal one. Our main goal is to find or estimate the values of the parameter $\varepsilon$, under which the spectrum of the operator $T(\varepsilon)$ is real and the operator itself is similar to a self-adjoint operator.

The main results of the paper are associated with the study of the operator $T(\varepsilon)$ with specific $\mathcal{PT}$-symmetric potential $p(x) = ix$. It turns out that this case presents  an exactly solvable model which allows us to trace the dynamics of the movement of the eigenvalues in all details and
to find explicitly the critical parameter values, in particular, to specify precisely the number $\varepsilon_1$ such that for $0<\varepsilon <\varepsilon_1$ the operator $T(\varepsilon)$ has a real spectrum and is similar to a self-adjoint one, but for  $\varepsilon \geq \varepsilon_1$ this property is  broken down.
In the general case, of course, to find explicitly critical
values is not possible. However, it is possible to get some  estimates for such  values of the parameter.

{\bf 2. General results.} The following statement shows that for sufficiently small values of the parameter $\varepsilon$
the operator $T(\varepsilon)$ is similar to a self-adjoint one.

\begin{Theorem}
\label{theorem1}
The spectrum of the operator $T(\varepsilon)$ consists of simple real eigenvalues, and the operator itself is similar to a self-adjoint one provided that any of the following conditions holds:
\begin{equation*}
\varepsilon < \frac{C_\infty}{\|p\|_\infty}, \quad
\varepsilon < \frac{C_1}{\|p\|_1}, \quad
\varepsilon < \frac{C_2}{\|p\|_2}.
\end{equation*}
Here $\|\cdot\|_\infty, \ \|\cdot\|_1, \ \|\cdot\|_2 $ are the norms in spaces $L_\infty(-1,1), \, L_1(-1,1)\,\, L_2(-1,1)$, respectively (assuming that the potential $p$ belongs to the corresponding space).
The constants in these inequalities are defined as follows
\begin{equation*}
C_\infty = \frac {3\pi^2}8;
\end{equation*}
\begin{equation*}
C_1=\frac{\pi^2}4 \left( \frac 23 + \sum_{k=2}^\infty
\left(k^2- \frac 52\right)^{-1}\right)^{-1};
\end{equation*}
\begin{equation*}
C_2=\frac{\pi^2}4 \left( \frac 49 +\sum_{k=2}^\infty
\left(k^2 - \frac 52\right)^{-2}\right)^{-1/2}.
\end{equation*}
\end{Theorem}

The proof of this theorem uses estimates for the resolvent of the operator family $L(\varepsilon)= \varepsilon T(\varepsilon)$, which can be considered as  a perturbation   of the operator $L(0)= -d^2/dx^2$ with boundary conditions \eqref{2}
(see \cite[$\S \, 7,9$]{Sh1}). It ia also used the Dunford theorem on the unconditional basis property of eigenfunctions of the operator
$T(\varepsilon)$ and theorems of Bari and Boas about the properties of unconditional bases \cite[$\S\, 6$]{Sh1}.

In the case of non-real  $\mathcal{PT}$-symmetric potential $p$ the  property of the  operator to be self-adjoint
 breaks down  as the parameter $\varepsilon$ is growing up.   For the case when the potential $p$ is a non-real polynomial the paper \cite{TS}
contains the result, which approves the localization of the spectrum along certain curves in the complex plane,
the structure of which is determined by the coefficients of the polynomial.  It follows from  this result that
  for large values of the parameter $\varepsilon$ the non-real eigenvalues do appear and their number
   increases proportionally to
 $\varepsilon^{1/2}$ as $\varepsilon\to\infty$.

{\bf 3. The complex Airy operator, as an explicitly solvable model.} Further we study the model operator
\begin{equation}
\label{3}
T(\varepsilon)\, y=-\varepsilon ^{-1} y''+ixy,\quad y(-1)=y(1)=0.
\end{equation}
 The  behavior of the spectrum of this operator  as $\varepsilon\to+\infty$ was investigated
in details   in the papers  \cite {Sh2, Sh3}. For large values of the parameter $\varepsilon$ the eigenvalues are concentrated along segments
$[i,1/\sqrt{3}]$, $[-i,1/\sqrt{3}]$ and the ray $[1/\sqrt{3},+\infty)$, and together they form the so-called {\it limit spectral graph}.
The point of $1/\sqrt{3}$ is  called the {\it knot point} of this limit spectral graph.
From theorem \ref{theorem1} we obtain that for the values
$\varepsilon < 3\pi^2/8 \approx 3.7011$
the spectrum of the operator \eqref{3} is real and consists of simple positive eigenvalues
$\lambda_k = \lambda_k(\varepsilon)$, which we assume to be  numerated in the increasing order. In this case
$\lambda_k (\varepsilon) \sim \varepsilon^{-1}(\pi k/2)^2 $
 as $\varepsilon \to 0$ and  $k$ is fixed, and the same is true as $k\to\infty$ and $\varepsilon$ is fixed.
  However, the estimate  $\varepsilon <\varepsilon_0 = 3\pi^2/8$
is too coarse. This estimate could be close to be sharp if the first two eigenvalues
$\mu_1(\varepsilon)$ and $\mu_2(\varepsilon)$ of the operator family  $L(\varepsilon) = \varepsilon T(\varepsilon)$
would move towards each other with increasing $\varepsilon \in [0,\, \varepsilon_0)$.  But it is  not the case.

In reality, the portrait  of the movements of the eigenvalues is the following. Numerical calculations show (in Theorem \ref{theorem4}
they get analytical confirmation) that
starting from small values of the parameter $\varepsilon$, all the eigenvalues of the
operator family $T(\varepsilon)$ move with increasing $\varepsilon >0$  from infinity to the left.
When $\varepsilon \approx 5.1$ the first eigenvalue crosses the knot point  $1/\sqrt 3\approx 0.58$
and continues to move to the left. When $\varepsilon =\varepsilon_{1,turn} \approx 9.3$
the first eigenvalue reaches some point $\lambda_{1, turn}\approx 0.45$,
stops at this point, and upon further increasing $\varepsilon >
\varepsilon_{1,turn}$ begins to move in the opposite direction, while all other eigenvalues continue to move to the left.
\begin{center}
\includegraphics[width=7cm,keepaspectratio]{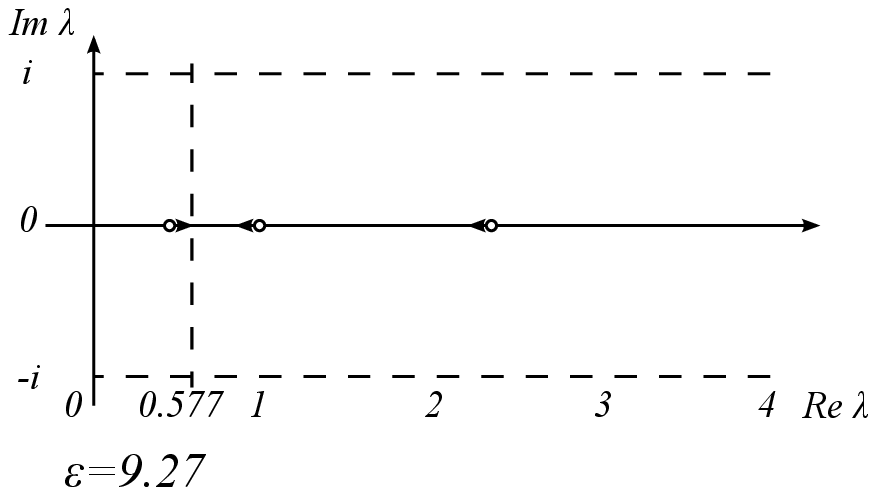}
\end{center}
Further, the first and second eigenvalues move towards each other approaching
 the knot point  $1/\sqrt 3$ and collide  at this point when $\varepsilon_1\approx12.3$.
\begin{center}
\includegraphics[width=7cm,keepaspectratio]{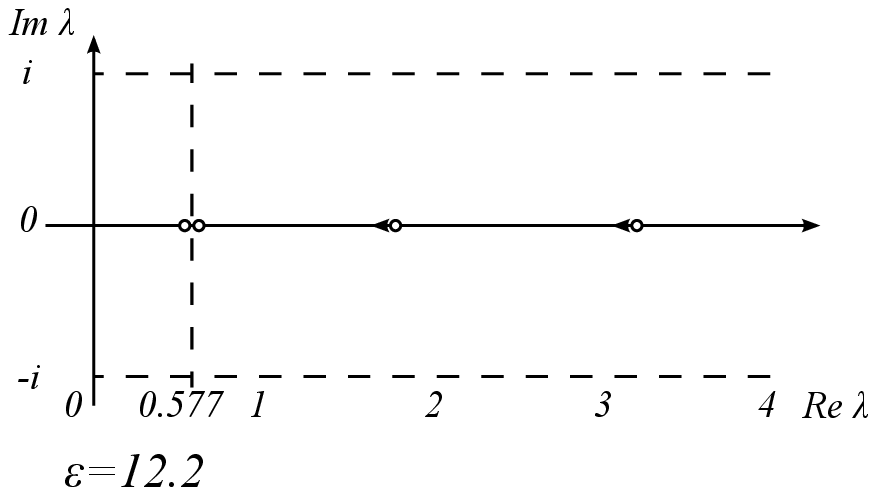}
\end{center}
After the collision in the knot point the first and second eigenvalues come off  in the complex plane
at the right angle to the real axis,
and upon further increasing $\varepsilon$, approach rapidly the segments $[1/\sqrt 3, \pm i]$ and continue the movement,
clinging to these segments, in the direction of the points $\pm i$.
\begin{center}
\includegraphics[width=7cm,keepaspectratio]{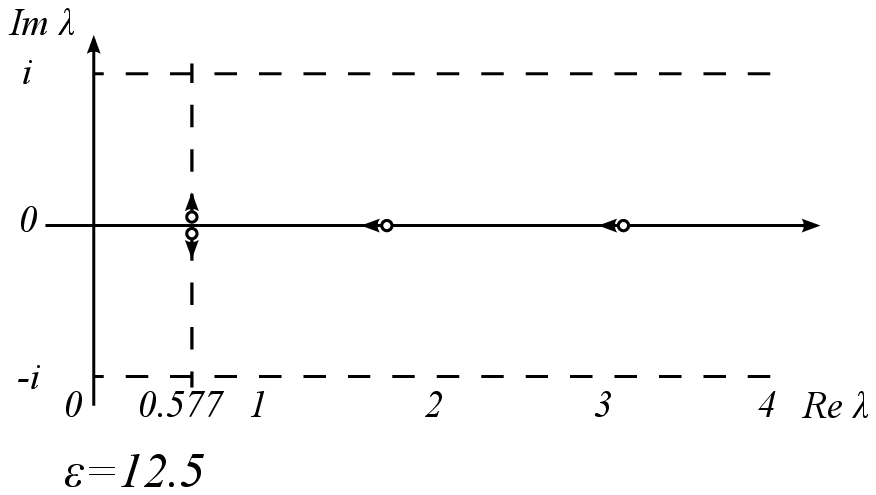}
\end{center}

Further, with growth of $\varepsilon$  the third eigenvalue crosses the knot point moving to the left to the point $\lambda_{2, turn} < 1/\sqrt 3$,
stopped at this point, and then moves in the opposite direction towards the fourth eigenvalue until the collision
 again at the knot point  $1/\sqrt 3$, and subsequently jumping in the complex plane. The fifth and sixth eigenvalues (and subsequent $2k-1$and $2k$th  ones)
repeat the same dynamics. For large $\varepsilon$ the eigenvalues accumulate on the real ray $[1/\sqrt 3, +\infty)$,
while all the non-real eigenvalues
nestle  to the segments $[1/\sqrt 3, \pm i]$, moving to the points $\pm i$.
\begin{center}
\includegraphics[width=7cm,keepaspectratio]{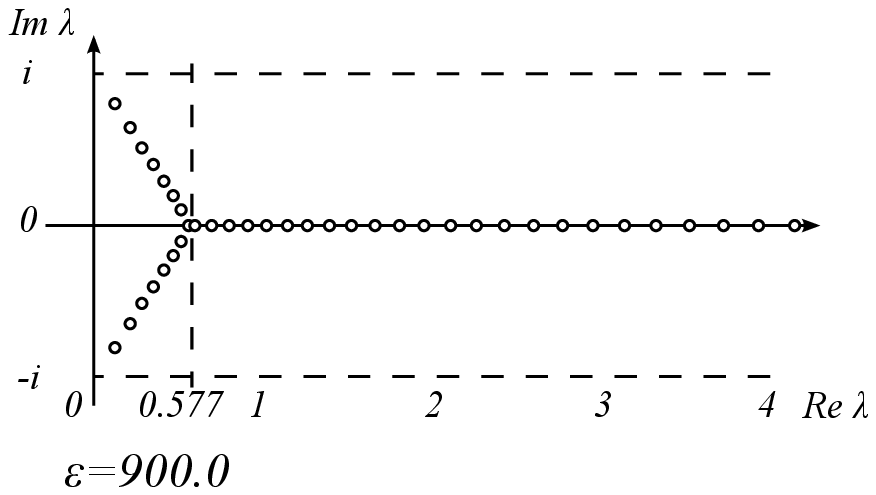}
\end{center}
Explicit formulae for the distribution of the eigenvalues in the intervals $[1/\sqrt 3, \pm i]$ and the ray
$[1/\sqrt 3, +\infty)$  for large parameter $\varepsilon$ can be  found in \cite{Sh2}.

%\begin{Definition}\n \ rm
Further, the set $\mathcal{E}$ of pairs $(\varepsilon,\lambda)\in\mathbb{R^{+}}\times\mathbb{C}$ for which for a given parameter $\varepsilon>0$ the number $\lambda$ belongs to spectrum of
the operator $T(\varepsilon)$ we call {\it the spectral locus} of the  operator family $T(\varepsilon)$. The subset $\mathcal{E}_\mathbb{R}\subset\mathcal{E}$, which
corresponds to the real eigenvalues $\lambda$, is called {\it the real spectral locus} of the family $T(\varepsilon)$.
%\end{Definition}

The following figure shows
the real spectral locus of the model operator (the figure reflects the real computer calculations).
\begin{center}
\includegraphics[width=7cm,keepaspectratio]{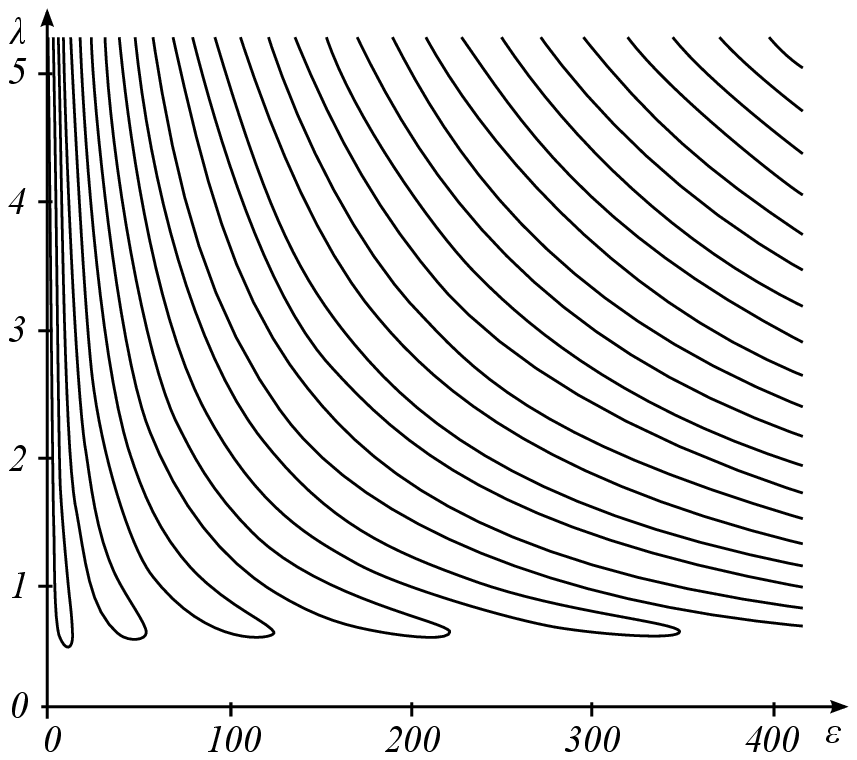}
\end{center}

\medskip
\medskip

Applying the Weierstrass preparation theorem to the characteristic determinant of the  eigenvalue problem for the operator $T(\varepsilon)$, we can prove the following statement.

\begin{Theorem}
\label{theorem2}
The real spectral locus of the family $T(\varepsilon)$ consists of regular analytic pairwise
non-intersecting Jordan  curves
in the extended complex plane with the ends in the infinity  point.
\end{Theorem}

A similar theorem was proved by Eremenko and Gabrielov~\cite{EG} while investigating the real spectrum of the cubic anharmonic oscillator on the whole axis.

From Theorem \ref{theorem2}, taking into account the smoothness of the curves  comprising the real part of the locus, and the absence of pairwise intersections, we get
the following statement.

\begin{Theorem}
\label{theorem3}
With increasing $\varepsilon$ the real eigenvalues of the operator family  \eqref{3} can collide only in pairs. Just before the moment of the collision the eigenvalues of a corresponding pair move towards each other, and immediately after the collision diverge in the complex plane at the right angles to the real axis.
\end{Theorem}

A point  $(\varepsilon_0,\lambda_0)\in\mathcal{E}$ of the spectral locus  we call {\it critical} if
the eigenvalue $\lambda_0$ of the operator $T(\varepsilon_0)$ is not simple (i.e. algebraic multiplicity of this
eigenvalue $> 1$).
Parameter values $\{\varepsilon_k\}_{k=1}^\infty$, for which at least one eigenvalue of the operator
$T(\varepsilon_k)$ is multiple,  we also call {\it critical} (it is easy to see that the set of such values is not more than countable).

{\bf 4. The main results for the model operator.}

Consider the classic Airy equation
$$
y'' =\xi\cdot y, \quad y =y(\xi),
$$
and its two standard solutions --- the functions $\Ai$ and $\Bi$
(see \cite{Olver}, for example).
A remarkable role in the sequel play the following special solutions of the Airy equation:
\begin{gather*}
U_-(\xi)=-\sqrt{3}\Ai(\xi)+\Bi(\xi),\\
U_+(\xi)=\sqrt{3}\Ai(\xi)+\Bi(\xi).
\end{gather*}

The following theorem is the main result of the present work.

\begin{Theorem}
\label{theorem4}
The zeros of the functions $U_-$ and $U_+$ are located on the rays $\arg z=\pi/3 + 2\pi k/3$, $k=-1,0,1$ symmetrically with respect to the origin. Let $\{\alpha_k\}_{k=1}^\infty$ and $\{\beta_k\}_{k=1}^\infty$ be the modules of the zeros of the functions $U_-$ and $U_+$, respectively, numerated  in the increasing order. These  zeros interlace:
$$
\alpha_0<\beta_1<\alpha_1<\beta_2<\alpha_2<\beta_3<\ldots.
$$
Denote
$$
\delta_k=\Bigl(\beta_k\frac{\sqrt{3}}{2}\Bigr)^3,\quad
\varepsilon_k=\Bigl(\alpha_k\frac{\sqrt{3}}{2}\Bigr)^3,\quad k\in\mathbb{N}
$$
(obviously, $0<\delta_1<\varepsilon_1<\delta_2<\varepsilon_2<\ldots$). Then

\begin{itemize}
\item  All the eigenvalues  of the operator $T(\varepsilon)$ are simple and this operator is similar to a self-adjoint
one, provided that $\varepsilon\in (0, \varepsilon_1)$, For all other values of $\varepsilon>0$ this property is
 broken down.
\item
All the critical points of the spectral locus coincide with the set
$$
\mathcal{M}=\Bigl\{
\left(\varepsilon_k,\frac{1}{\sqrt{3}}\right)
\Bigr\}.
$$
For all critical values of the parameter the knot point  $1/\sqrt{3}$ is 2-multiple eigenvalue of
 the operator $T(\varepsilon_k)$, which meets the Jordan cell.
 For all $\varepsilon \ne \varepsilon_k$, the operator $T(\varepsilon)$ is similar to a normal one  (i.e. to an  operator, commuting with its adjoint).
\item
All the odd eigenvalues $\lambda_{2k-1}(\varepsilon)$ move to the left (decrease) with increasing
parameter $\varepsilon$,  pass through  the knot  point $1/\sqrt 3$ as the parameter takes the values 
 $ \varepsilon = \delta_k$,
and continue to move to the left until they reach some points $\lambda_{2k-1,\, turn} < 1/\sqrt 3$.  Then, starting from these points,  they  turn back
and move to the right before colliding with the even  eigenvalues
$\lambda_{2k}$ in the knot  point $1/\sqrt 3$ as the parameter takes  the critical values  $\varepsilon = \varepsilon_k$.
\item
After the collision the eigenvalues move  in the opposite directions to the complex plane perpendicular to the real axis, and subsequently never come back to the real axis.
 Outside of the real axis the eigenvalues are unable to face.
\item
For the values  $\varepsilon=\delta_k$ and  $\lambda_{2k-1}=1/\sqrt{3}$ the eigenfunctions can be written
down explicitly:
$$
y(z)=U_+\left(\delta_k^{1/3}\left(\frac{1}{\sqrt{3}}-iz\right)\right).
$$
\item For the values  $\varepsilon=\varepsilon_k$ and  $\lambda_{2k-1}=\lambda_{2k}=1/\sqrt{3}$ the
eigenfunctions are also  written down explicitly:
$$
y(z)=U_-\left(\varepsilon^{1/3}\left(\frac{1}{\sqrt{3}}-iz\right)\right).
$$
\item
The following  asymptotics are valid  as $k\to\infty$:
\begin{gather*}
\varepsilon_k=\frac{\sqrt{3}}{4}\left(\frac{3}{2}\right)^3\left(\pi k-\frac{\pi}{12}+O\left(\frac{1}{k}\right)\right)^2,\\
\delta_k=\frac{\sqrt{3}}{4}\left(\frac{3}{2}\right)^3\left(\pi k-\frac{5\pi}{12}+O\left(\frac{1}{k}\right)\right)^2.
\end{gather*}

\end{itemize}
\end{Theorem}

The following theorem gives an upper bound for the turning points of the eigenvalues $\lambda_{2k-1, turn}$, $k\ge1$.

\begin{Theorem}
\label{theorem5}
Let the modules of the complex zeros $\{z_k\}_{k=1}^\infty$ of the function $\Bi$ lying in the first quadrant of the complex plane be numerated in the increasing order of their modulas. Then the following estimate  for the turning points of the odd eigenvalues is valid:
$$
\lambda_{2k-1, turn}<\cot\arg z_k<1/\sqrt{3}.
$$
\end{Theorem}

From this statement, in particular, we have $\lambda_{1,turn}<0.457$.

Let us say some words about the idea of the proof of Theorems \ref{theorem4} and \ref{theorem5}. It turns out that the 
substitution $\xi=\varepsilon^{1/3}(\lambda-ix)$ reduces the eigenvalue problem for the operator $T(\varepsilon)$ to the study
 of the zero distribution  of solutions of the Airy equation
$$
y'' = \xi y.
$$

The function $V_a(\xi)=a\Ai(\xi)+\Bi(\xi)$ of the argument $\xi\in\mathbb{C}$, with $a\in\mathbb{C}$, up to a constant multiplier, describes
all the solutions to this equation with the exception of $\Ai$.

If $a$ is real, then the complex (non-real) zeros of $V_a(\xi)$ specify an alternative parametrization of the real
 spectral locus of the original problem. Each point $(\varepsilon,\lambda)$ of this locus corresponds in a unique way
  to the pair of complex conjugate zeros
 of a real solutions of the Airy equation. Vice versa, if $\xi_0$ is zero of $V_a(\xi)$ for some real $a$, then
the real point of the spectral locus:
$\varepsilon=|\im\xi_0|^3$, $\lambda=\re\xi_0/|\im\xi_0|$ is uniquely determined. The equation $V_a(\xi)=0$ defines a countable number of implicit
functions  $\xi_k=\xi_k(a)$, $k\in\mathbb{Z}$, each of which admits an analytic  continuation in a neighborhood of the real axis and takes the values in the first or
fourth quadrants of the complex plane. The images  $\Gamma_k$ of the real axis of each of these functions
are pairwise disjoint analytic Jordan arcs,
unlimited with both ends.
\begin{center}
\includegraphics[width=7cm,keepaspectratio]{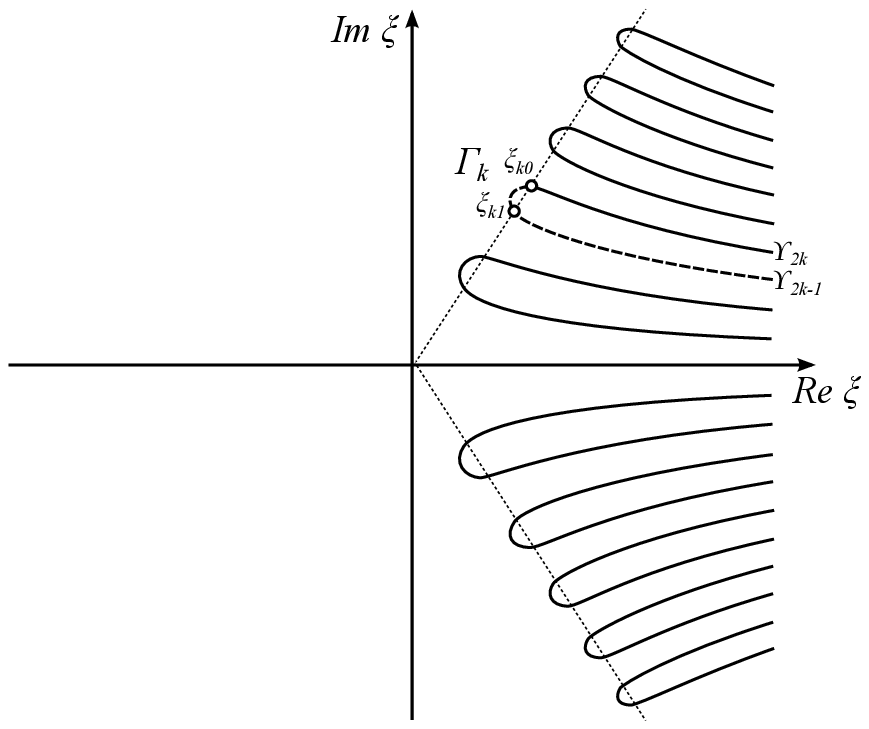}
\end{center}

\medskip
\medskip

Each of the ends of $\Gamma_k$, tending to infinity, approaches the real positive half-line.

If $\xi_k(a)$, $a\in\mathbb{R}$, is the parametrization of a curve $\Gamma_k$,
then the point $\xi_{k0}=\xi_k(a_0)$ of the curve $\Gamma_k$   corresponds to the critical point of the spectral locus of the operator family $T(\varepsilon)$ if and only if  $\xi_k'(a_0)\in\mathbb{R}$. This is true if and only if the tangent to the curve $\Gamma_k$ at the point $\xi_{k0}$ is parallel to the real axis. This is true if and only if  $a=-\sqrt{3}$ for all $k$.

For each $\Gamma_k$ lying in the first quadrant of the complex plane,
there exists the only  point $\xi_{k0}=\xi_k(-\sqrt{3})$  and it coincides with the zero of the functions $U_-$, lying
on the ray $l=\{\arg z=\pi/3\}$. Among  all the points of the curve  $\Gamma_k$, it has the largest imaginary part.

The point $\xi_{k0}$ divides $\Gamma_k$ on two arcs: $\gamma_{2k-1}$ and $\gamma_{2k}$, which correspond to the dynamics of the  eigenvalues $\lambda_{2k-1}$ and
$\lambda_{2k}$ of the original problem. The point $\xi_{k1}=\xi_k(\sqrt{3})$ lying on the ray $\arg z=\pi/3$ corresponds to passing of the $\lambda_{2k-1}$ through the knot point,
and the  $\xi_{k0}$ corresponds to collision of the $\lambda_{2k-1}$ and $\lambda_{2k}$.

The statements of the theorem concerning  the behavior of complex eigenvalues (the impossibility of coming back to the real axis and the impossibility of collisions outside
the real axis) are proved while investigating the characteristic determinant of the original problem in terms of the new variable $\xi$ and the parameter $\varepsilon$.

The work was supported by the Russian science Foundation, grant No. 17-11-01215. It will be published in Doklady Mathematics in 2017.

\end{document}